\documentclass[a4paper,11pt]{article}
\usepackage{amsmath, amssymb, amsthm}
\usepackage{graphicx} 

\begin{document}

\def\endproof{\hfill$\square$}
\def\PP{{\mathcal P}} \def\HH{{\mathcal H}} \def\MM{{\mathcal M}}
\def\D{{\mathbb D}} \def\C{{\mathbb C}} \def\Sym{{\rm Sym}}
\def\gr{{\rm gr}} \def\Z{{\mathbb Z}} \def\Id{{\rm Id}}
\def\S{{\mathbb S}} \def\T{{\mathbb T}} \def\N{{\mathbb N}}
\def\R{{\mathbb R}} \def\St{{\rm St}} \def\pos{{\rm pos}}
\def\An{{\rm An}} \def\OO{{\mathcal O}} \def\DD{{\mathcal D}}
\def\SS{{\mathcal S}} \def\PP{{\mathcal P}}\def\AA{{\mathcal A}}
\def\CC{{\mathcal C}} \def\A{{\mathbb A}}

\title{\bf{From braid groups to mapping class groups}}
 
\author{\textsc{Luis Paris}}

\date{\today}

\maketitle

\begin{abstract}
\noindent
This paper is a survey of some properties of the braid groups and related groups that 
lead to questions on mapping class groups.
\end{abstract}


\section{Introduction}

I have to confess that I am not a leader in the area of mapping class groups and 
Teichm\"uller geometry. Nevertheless, I can ask a single question which, I believe, should be of 
interest for experts in the field.

\bigskip\noindent
{\bf Question.} {\it Which properties of the braid groups can be extended to the mapping class 
groups?}

\bigskip\noindent
The present paper is a sort of survey of some (not always well-known) properties of the braid groups 
and related groups which lead to questions on mapping class groups. The reader will find 
throughout the text several more or less explicit questions, but none of them are stated as 
``official'' problems or conjectures. This choice is based on the idea that the proof often comes 
together with the question (or the theorem), and, moreover, the young researchers should keep in mind 
that the surroundings and the applications of a given problem should be the main motivations for 
its study, and not the question itself (or its author).

\bigskip\noindent
Throughout the paper, $\Sigma$ will denote a compact, connected, and oriented surface. The boundary 
$\partial \Sigma$, if non empty, is a finite collection of simple closed curves. Consider a 
finite set $\PP=\{P_1, \dots, P_n\}$ of $n$ distinct points, called {\it punctures}, in the 
interior of $\Sigma$. Define $\HH( \Sigma, \PP)$ to be the group of orientation-preserving 
homeomorphisms $h: \Sigma \to \Sigma$ such that $h$ is the identity on each boundary component of 
$\Sigma$ and $h(\PP)=\PP$. The {\it mapping class 
group $\MM(\Sigma, \PP)= \pi_0( \HH( \Sigma, \PP))$ of $\Sigma$ relative to $\PP$} is defined to be the set of 
isotopy classes of mappings in $\HH( \Sigma, \PP)$, with composition as group operation. We 
emphasize that, throughout an isotopy, the boundary components and the punctures of $\PP$ remain 
fixed. It is clear that, up to isomorphism, this group depends only on the genus $g$ of $\Sigma$, 
on the number $p$ of components of $\partial \Sigma$, and on the cardinality $n$ of $\PP$. So, 
we may write $\MM( g,p,n)$ in place of $\MM( \Sigma, \PP)$, and simply $\MM(\Sigma)$ for $\MM( 
\Sigma, \emptyset)$.

\bigskip\noindent
Let $\D=\{z \in \C; |z|\le 1\}$ be the {\it standard disk}, and let $\PP=\{P_1, \dots, P_n\}$ be 
a finite collection of $n$ punctures in the interior of $\D$. Define a {\it $n$-braid} to be a $n$-tuple 
$\beta=( b_1, \dots, b_n)$ of disjoint smooth paths in $\D \times [0,1]$, called the {\it 
strings} of $\beta$, such that:
\begin{itemize}
\item the projection of $b_i(t)$ on the second coordinate is $t$, for all $t \in [0,1]$ and all 
$i \in \{1, \dots, n\}$;
\item $b_i(0)=(P_i,0)$ and $b_i(1)= (P_{\chi(i)},1)$, where $\chi$ is a permutation of $\{1, \dots, 
n\}$, for all $i \in \{1, \dots, n\}$.
\end{itemize}
An isotopy in this context is a deformation through braids which fixes the ends. Multiplication of 
braids is defined by concatenation. The isotopy classes of braids with this multiplication form a 
group, called the {\it braid group on $n$ strings}, and denoted by $B_n(\PP)$. This group 
does not depend, up to isomorphism, on the set $\PP$, but only on the cardinality $n$ of $\PP$, 
thus we may write $B_n$ in place of $B_n(\PP)$.

\bigskip\noindent
The group $B_n$ has a well-known presentation with generators $\sigma_1, \dots, \sigma_{n-1}$ 
and relations
\[
\begin{array}{cl}
\sigma_i \sigma_j = \sigma_j \sigma_i &\quad \text{if } |j-i|>1\,,\\
\sigma_i \sigma_j \sigma_i = \sigma_j \sigma_i \sigma_j &\quad\text{if } |j-i|=1\,.
\end{array}
\]
There are other equivalent descriptions of $B_n$, as a group of automorphisms of the 
free group $F_n$, as the fundamental group of a configuration space, or as the mapping class group
of the $n$-punctured disk. 
This explains the importance of the braid groups in many disciplines.

\bigskip\noindent
Now, the equality
\begin{equation}\label{eq11}
B_n(\PP)=\MM(\D,\PP)
\end{equation}
is the source of the several questions that the reader will find in the paper.


\section{Garside groups and Artin groups}

We start this section with a brief presentation of Garside groups. Our presentation of the 
subject draws in many ways from the work of Dehornoy \cite{Deh1} as well as \cite{DePa}, and, 
like all treatments of Garside groups, is inspired ultimately by the seminal papers of Garside 
\cite{Gar}, on braid groups, and Brieskorn and Saito \cite{BrSa}, on Artin groups. 

\bigskip\noindent
Let $M$ be an arbitrary monoid. We say that $M$ is {\it atomic} if there exists a function $\nu: 
M \to \N$ such that
\begin{itemize}
\item $\nu(a)=0$ if and only if $a=1$;
\item $\nu(ab) \ge \nu(a) + \nu(b)$ for all $a,b \in M$.
\end{itemize}
Such a function on $M$ is called a {\it norm} on $M$.

\bigskip\noindent
An element $a \in M$ is called an {\it atom} if it is indecomposable, namely, if $a=bc$ then either 
$b=1$ or $c=1$. This definition of atomicity is taken from \cite{DePa}. We refer to \cite{DePa} 
for a list of further properties all equivalent to atomicity. In the same paper, it is shown that 
a subset of $M$ generates $M$ if and only if it contains the set of all atoms. In particular, 
$M$ is finitely generated if and only if it has finitely many atoms.

\bigskip\noindent
Given that a monoid $M$ is atomic, we may define left and right invariant partial orders $\le_L$ 
and $\le_R$ on $M$ as follows.
\begin{itemize}
\item Set $a \le_L b$ if there exists $c \in M$ such that $ac=b$.
\item Set $a \le_R b$ if there exists $c \in M$ such that $ca=b$.
\end{itemize}
We shall call these {\it left} and {\it right divisibility orders} on $M$.

\bigskip\noindent
A {\it Garside monoid} is a monoid $M$ such that
\begin{itemize}
\item $M$ is atomic and finitely generated;
\item $M$ is cancellative;
\item $(M, \le_L)$ and $(M, \le_R)$ are lattices;
\item there exists an element $\Delta \in M$, called {\it Garside element}, such that the set 
$L(\Delta)= \{ x \in M; x\le_L \Delta\}$ generates $M$ and is equal to $R(\Delta)= \{ x\in M; 
x\le_R \Delta\}$.
\end{itemize}
For any monoid $M$, we can define the group $G(M)$ which is presented by the generating set $M$ 
and the relations $ab=c$ whenever $ab=c$ in $M$. There is an obvious canonical homomorphism $M 
\to G(M)$. This homomorphism is not injective in general. The group $G(M)$ is known as the {\it 
group of fractions} of $M$. Define a {\it Garside group} to be the group of fractions of a 
Garside monoid.

\bigskip\noindent
{\bf Remarks.}
\begin{enumerate}
\item 
Recall that a monoid $M$ satisfies {\it \"Ore's conditions} if $M$ is cancellative and if,
for all $a,b \in M$, there exist $c,d \in M$ such that $ac=bd$. If $M$ satisfies \"Ore's conditions
then the canonical homomorphism $M \to G(M)$ is injective. A Garside monoid obviously
satisfies \"Ore's conditions, thus any Garside monoid embeds in its group of fractions.
\item 
A Garside element is never unique. For example, if $\Delta$ is a Garside element, then 
$\Delta^k$ is also a Garside element for all $k \ge 1$ (see \cite{Deh1}).
\end{enumerate}
Let $M$ be a Garside monoid. The lattice operations of $(M,\le_L)$ are denoted by $\vee_L$ and 
$\wedge_L$. 
These are defined as follows. For $a,b,c \in M$, $a\wedge_L b \le_L a$, $a \wedge_L b \le_L b$, and
if $c\le_L a$ and $c \le_L b$, then $c \le_L a \wedge_L b$. For $a,b,c \in M$, $a\le_L a\vee_L b$,
$b \le_L a \vee_L b$ and, if $a \le_L c$ and $b \le_L c$, then $a \vee_L b \le_L c$.
Similarly, the lattice operations of $(M,\le_R)$ are denoted by $\vee_R$ and 
$\wedge_R$.

\bigskip\noindent
Now, we briefly explain how to define a biautomatic structure on a given Garside group. By 
\cite{EpAl}, such a structure furnishes solutions to the word problem and to the conjugacy 
problem, and it implies that the group has quadratic isoperimetric inequalities. We refer to 
\cite{EpAl} for definitions and properties of automatic groups, and to \cite{Deh1} for more 
details on the biautomatic structures on Garside groups.

\bigskip\noindent
Let $M$ be a Garside monoid, and let $\Delta$ be a Garside element of $M$. For $a \in M$, we 
write $\pi_L(a) = \Delta \wedge_L a$ and denote by $\partial_L(a)$ the unique element of $M$ such 
that $a=\pi_L(a) \cdot \partial_L(a)$. Using the fact that $M$ is atomic and that $L(\Delta)= 
\{x \in M; x \le_L \Delta\}$ generates $M$, one can easily show that $\pi_L(a) \neq 1$ if $a \neq 
1$, and that there exists a positive integer $k$ such that $\partial_L^k(a)=1$. Let $k$ be the 
lowest integer satisfying $\partial_L^k(a)=1$. Then the expression
\[
a=\pi_L(a) \cdot \pi_L(\partial_L(a)) \cdot \dots \cdot \pi_L(\partial_L^{k-1}(a))
\]
is called the {\it normal form} of $a$.

\bigskip\noindent
Let $G=G(M)$ be the group of fractions of $M$. Let $c \in G$. Since $G$ is a lattice with 
positive cone $M$, the element $c$ can be written $c=a^{-1}b$ with $a,b \in M$. Obviously, $a$ 
and $b$ can be chosen so that $a \wedge_L b =1$ and, with this extra condition, are unique. Let 
$a=a_1a_2 \dots a_p$ and $b=b_1b_2 \dots b_q$ be the normal forms of $a$ and $b$, respectively. 
Then the expression
\[
c=a_p^{-1} \dots a_2^{-1} a_1^{-1} b_1b_2 \dots b_q
\]
is called the {\it normal form} of $c$.

\bigskip\noindent
{\bf Theorem 2.1 (Dehornoy \cite{Deh1}).} {\it Let $M$ be a Garside monoid, 
and let $G$ be the group of fractions of $M$. Then the normal forms of the elements of $G$ form 
a symmetric rational language on the finite set $L(\Delta)$ which has the fellow traveler 
property. In particular, $G$ is biautomatic.}

\bigskip\noindent
The notion of a Garside group was introduced by Dehornoy and the author \cite{DePa} in a slightly 
restricted sense, and, later, by Dehornoy \cite{Deh1} in the larger sense which is now generally 
used. As pointed out before, the theory of Garside groups is largely inspired by the papers of 
Garside \cite{Gar}, which treated the case of the braid groups, and Brieskorn and Saito 
\cite{BrSa}, which generalized Garside's work to Artin groups. The Artin groups of spherical 
type which include, notably, the braid groups, are motivating examples. Other interesting 
examples include all torus link groups (see \cite{Pic1}) and some generalized braid groups 
associated to finite complex reflection groups (see \cite{BesCor}).

\bigskip\noindent
Garside groups have many attractive properties: solutions to the word and conjugacy problems 
are known (see \cite{Deh1}, \cite{Pic2}, \cite{FrGo}), they are torsion-free (see 
\cite{Deh2}), they are biautomatic (see \cite{Deh1}), and they admit finite dimensional 
classifying spaces (see \cite{DeLa}, \cite{CMW}). 
There also exist criteria in terms of presentations which detect Garside groups 
(see \cite{DePa}, \cite{Deh1}).

\bigskip\noindent
Let $S$ be a finite set. A {\it Coxeter matrix} over $S$ is a square matrix 
$M=(m_{s\,t})_{s,t\in S}$ indexed by the elements of $S$ and such that $m_{s\,s}=1$ for all $s 
\in S$, and $m_{s\,t}=m_{t\,s} \in \{2,3,4, \dots, +\infty\}$ for all $s,t \in S$, $s \neq t$. 
A Coxeter matrix $M=(m_{s\,t})_{s,t\in S}$ is usually represented by its {\it Coxeter graph}, 
$\Gamma$, which is defined as follows. The set of vertices of $\Gamma$ is $S$, two vertices $s,t$ 
are joined by an edge if $m_{s\,t} \ge 3$, and this edge is labeled by $m_{s\,t}$ if $m_{s\,t} 
\ge 4$.

\bigskip\noindent
Let $\Gamma$ be a Coxeter graph with set of vertices $S$. For two objects $a,b$ and $m \in \N$, 
define the word
\[\omega(a,b:m)= \left\{
\begin{array}{ll}
(ab)^{{m \over 2}}\quad&\text{if }m\text{ is even}\,,\\
(ab)^{{m-1\over 2}}a \quad &\text{if }m\text{ is odd}\,.
\end{array}\right.
\]
Take an abstract set $\SS= \{ \sigma_s; s \in S\}$ in one-to-one correspondence with $S$. The {\it 
Artin group of type $\Gamma$} is the group $A=A_\Gamma$ generated by $\SS$ and subject to the 
relations
\[
\omega(\sigma_s, \sigma_t: m_{s\,t}) = \omega( \sigma_t, \sigma_s: m_{s\,t}) \quad \text{for }s,t 
\in S,\ s \neq t,\text{ and } m_{s\,t}< +\infty\,,
\]
where $M=(m_{s\,t})_{s,t\in S}$ is the Coxeter matrix of $\Gamma$. The {\it Coxeter group of type 
$\Gamma$} is the group $W=W_\Gamma$ generated by $S$ and subject to the relations
\begin{gather*}
s^2=1 \quad \text{for } s\in S\,,\\
(st)^{m_{s\,t}}=1 \quad \text{for }s,t\in S,\ s\neq t,\text{ and }m_{s\,t}<+\infty\,.
\end{gather*}
Note that $W$ is the quotient of $A$ by the relations $\sigma_s^2=1$, $s \in S$, and $s$ is the 
image of $\sigma_s$ under the quotient epimorphism.

\bigskip\noindent
The number $n=|S|$ of generators is called the {\it rank} of the Artin group (and of the Coxeter 
group). We say that $A$ is {\it irreducible} if $\Gamma$ is connected, and that $A$ is of 
{\it spherical type} if $W$ is finite.

\bigskip\noindent
Coxeter groups have been widely studied. Basic references for them are \cite{Bou} and 
\cite{Hum}. In contrast,
Artin groups are poorly understood in general. Beside the spherical type ones, 
which are Garside groups,
the Artin groups 
which are better understood include right--angled Artin groups, 2-dimensional Artin groups, and 
FC type Artin groups. Right-angled Artin groups (also known as graph groups or free partially 
commutative groups) have been widely studied, and their applications extend to various 
domains such as parallel computation, random walks, and cohomology of groups. 2-dimensional Artin 
groups and FC type Artin groups have been introduced by Charney and Davis \cite{ChDa} in 1995 in 
their study of the $K(\pi,1)$-problem for complements of infinite hyperplane arrangements 
associated to reflection groups.

\bigskip\noindent
Let $G$ be an Artin group (resp. Garside group). Define a {\it geometric representation} of $G$ 
to be a homomorphism from $G$ to some mapping class group. 
Although I do not believe that mapping class groups are either
Artin groups or Garside groups, the theory of geometric representations
should be of interest from the perspective of both
mapping class groups and Artin (resp. Garside) groups.
For instance, it is not known which Artin groups (resp. Garside groups) can be embedded into a
mapping class group.

\bigskip\noindent
A first example of such representations are the (maximal) free abelian subgroups of the mapping class groups
whose study was 
initiated by Birman, Lubotzky, and McCarthy \cite{BLM}. These representations are crucial in the study of the 
automomorphism groups and abstract commensurators of the mapping class groups (see \cite{Iva1}, 
\cite{Iva2}, \cite{IvMc}, \cite{McM2}). Another example are the presentations of the mapping class groups as 
quotients of Artin groups by some (not too bad) relations (see \cite{Mat}, \cite{LaPa}). 
However, the most attractive 
examples of geometric representations are the so-called {\it monodromy 
representations} introduced by several people (Sullivan, Arnold, A'Campo) for the Artin groups of 
type $A_n$, $D_n$, $E_6$, $E_7$, $E_8$, and extended by Perron, Vannier \cite{PeVa}, Crisp 
and the author \cite{CrPa} to all small type Artin groups as follows.

\bigskip\noindent
First, recall that an Artin group $A$ associated to a Coxeter graph $\Gamma$ 
is said to be of {\it small type} if $m_{s\,t} \in \{2,3\}$ for all $s,t \in S$, $s \neq t$,
where $M=(m_{s\,t})_{s,t\in S}$ is the Coxeter matrix of $\Gamma$.

\bigskip\noindent
An {\it essential circle} in $\Sigma$ is an oriented embedding $a: \S^1 \to \Sigma$ such that 
$a(\S^1) \cap \partial \Sigma = \emptyset$, and $a(\S^1)$ does not bound any disk. We shall use 
the description of the circle $\S^1$ as $\R/2\pi \Z$. Take an (oriented) embedding $\A_a: \S^1 
\times [0,1] \to \Sigma$ of the annulus such that $\A_a(\theta, {1 \over 2}) = a(\theta)$ for 
all $\theta \in \S^1$, and define a homomorphism $T_a \in \HH (\Sigma)$, which restricts to the 
identity outside the interior of the image of $\A_a$, and is given by
\[
(T_a \circ \A_a) (\theta, x)= \A_a( \theta +2 \pi x,x) \quad \text{for } (\theta,x) \in \S^1 
\times [0,1]\,,
\]
inside the image of $\A_a$.
The {\it Dehn twist} along $a$ is the element of $\MM(\Sigma)$ represented by $T_a$. The following 
result is easily checked and may be found, for instance, in \cite{Bir2}.

\bigskip\noindent
{\bf Proposition 2.2.} {\it Let $a_1, a_2: \S^1 \to \Sigma$ be two essential circles, and, for 
$i=1,2$, let $\tau_i$ denote the Dehn twist along $a_i$. Then
\[
\begin{array}{cl}
\tau_1 \tau_2 = \tau_2 \tau_1 &\quad \text{if }a_1\cap a_2 = \emptyset\,,\\
\tau_1 \tau_2 \tau_1 = \tau_2 \tau_1 \tau_2 &\quad \text{if } |a_1 \cap a_2|=1\,.
\end{array}
\]}

\noindent
We assume now that $\Gamma$ is a small type Coxeter graph, and we 
associate to $\Gamma$ a surface $\Sigma= \Sigma(\Gamma)$ as follows.

\bigskip\noindent
Let $M=(m_{s\,t})_{s,t\in S}$ be the Coxeter matrix of $\Gamma$.
Let $<$ be a total order on $S$ which can be chosen arbitrarily. For each $s \in S$, we define 
the set
\[
\St_s= \{t \in S; m_{s\,t}=3\} \cup \{s\}\,.
\]
We write $\St_s= \{t_1, t_2, \dots, t_k\}$ such that $t_1<t_2< \dots <t_k$, and suppose $s=t_j$. 
The difference $i-j$ is called the {\it relative position} of $t_i$ with respect to $s$, and is 
denoted by $\pos(t_i:s)$. In particular, $\pos(s:s)=0$.

\bigskip\noindent
Let $s\in S$. Put $k=|St_s|$. We denote by $\An_s$ the annulus defined by
\[
\An_s= \R /2k\Z \times [0,1]\,.
\]
For each $s \in S$, write $P_s$ for the point $(0,0)$ of $\An_s$. The surface $\Sigma =\Sigma(\Gamma)$ 
is defined by
\[
\Sigma= \left( \coprod_{s \in S} \An_s \right) / \approx\,,
\]
where $\approx$ is the relation defined as follows. Let $s,t \in S$ such that $m_{s\,t}=3$ and 
$s<t$. Put $p=\pos(t:s)>0$ and $q=\pos(s:t)<0$. For each $(x,y) \in [0,1] \times [0,1]$, the 
relation $\approx$ identifies the point $(2p+x,y)$ of $\An_s$ with the point $(2q+1-y,x)$ of 
$\An_t$ (see Figure 1). We identify each annulus $\An_s$ and the point $P_s$ with their image in 
$\Sigma$, respectively.

\begin{figure}
\begin{center}
\includegraphics[width=14cm]{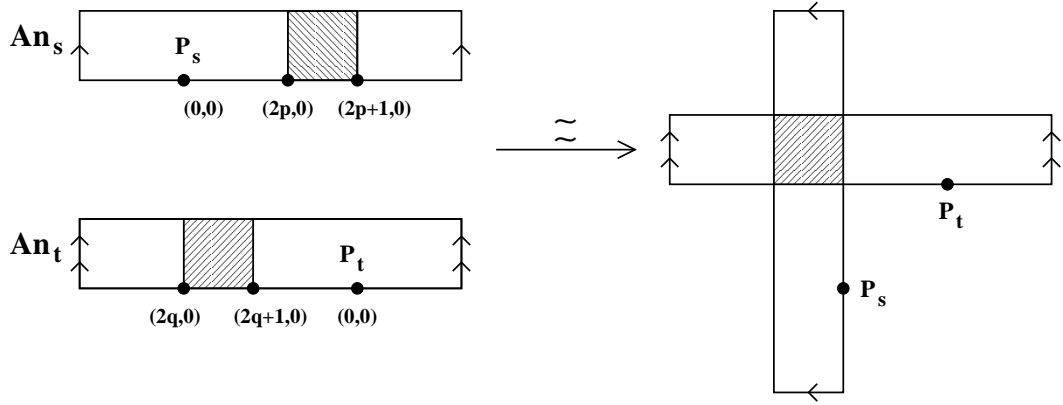}
\end{center}
\caption{Identification of annuli.}
\end{figure}

\bigskip\noindent
We now define the monodromy representation $\rho: A_\Gamma \to \MM (\Sigma)$. Let $s \in S$, and put 
$k=|\St_s|$. We denote by $a_s: \S^1 \to \Sigma$ the essential circle of $\Sigma$ such that 
$a_s(\theta)$ is the point $({k \theta \over \pi}, {1 \over 2})$ of $\An_s$ (see Figure 2). We 
let $\tau_s$ denote the Dehn twist along $a_s$. One has $a_s \cap a_t = \emptyset$ if 
$m_{s\,t}=2$, and $|a_s \cap a_t|=1$ if $m_{s\,t}=3$. Therefore, by Proposition 2.2, we have the 
following.

\begin{figure}
\begin{center}
\includegraphics[width=7cm]{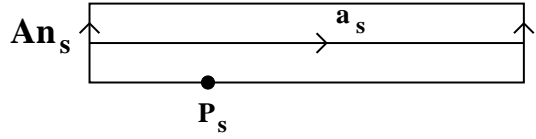}
\end{center}
\caption{The essential circle $a_s$ in $\An_s$.}
\end{figure}

\bigskip\noindent
{\bf Proposition 2.3.} {\it Let $\Gamma$ be a small type Coxeter graph. There exists a well-defined 
group homomorphism $\rho: A_\Gamma \to \MM (\Sigma)$ which sends $\sigma_s$ to $\tau_s$ for each $s 
\in S$.}

\bigskip\noindent
{\bf Example.} Consider the braid group $B_n$, $n \ge 3$. If $n$ is odd, then 
the associated surface, $\Sigma$, is a surface of genus ${n-1 \over 2}$ with one boundary 
component. If $n$ is even, then $\Sigma$ is a surface of genus ${n-2 \over 2}$ with two 
boundary components (see Figure 3). For each $i=1, \dots, n-1$, let $a_i: \S^1 \to \Sigma$ denote 
the essential circle pictured in Figure 3, and let $\tau_i$ denote the Dehn twist along $a_i$. 
Then the monodromy representation $\rho: B_n \to \MM (\Sigma)$ of the braid group sends 
$\sigma_i$ to $\tau_i$ for all $i=1, \dots, n-1$.

\begin{figure}
\begin{center}
\includegraphics[width=11cm]{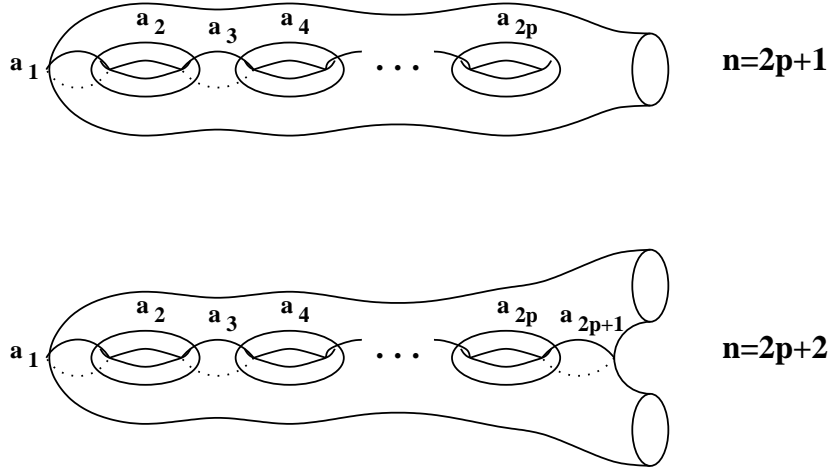}
\end{center}
\caption{Monodromy representation for $B_n$.}
\end{figure}

\bigskip\noindent
Let $\OO_2$ denote the ring of germs of functions at $0$ of $\C^2$.
Let $f \in \OO_2$ with an isolated singularity at $0$, and let $\mu$ be the Milnor number of $f$. 
To the germ $f$ one can associate a surface $\Sigma$ with boundary, called the {\it Milnor fiber} 
of $f$, an analytic subvariety $\DD$ of $\D (\eta)= \{t \in \C^\mu; \| t\| <\eta \}$, where 
$\eta >0$ is a small number, and a representation $\rho: \pi_1( \D(\eta) \setminus \DD) \to \MM 
(\Sigma)$, called the {\it geometric monodromy} of $f$. Recall that the simple singularities are 
defined by the equations
\[
\begin{array}{ll}
A_n: \quad&f(x,y)=x^2+y^{n+1}\,, \quad n\ge 1\,,\\
D_n: \quad&f(x,y)=x(x^{n-2}+y^2)\,, \quad n\ge 4\,,\\
E_6:\quad&f(x,y)=x^3+y^4\,,\\
E_7: \quad&f(x,y)=x(x^2+y^3)\,,\\
E_8:\quad&f(x,y)=x^3+y^5\,.
\end{array}
\]
By \cite{Arn} and \cite{Bri}, for each $\Gamma \in \{A_n; n\ge 1\} \cup \{D_n; n\ge 4\} \cup 
\{E_6,E_7,E_8\}$, the group $\pi_1( \D(\eta) \setminus \DD)$ is isomorphic to the Artin group of 
type $\Gamma$, and $\rho: A_\Gamma \to \MM( \Sigma)$ is the monodromy representation of the group $A_\Gamma$.

\bigskip\noindent
The monodromy representations are known to be faithful for $\Gamma=A_n$ and $\Gamma=D_n$ (see 
\cite{PeVa}), but are not faithful for $\Gamma \in \{E_6,E_7,E_8\}$ (see \cite{Waj}). Actually, the monodromy 
representations are not faithful in general (see \cite{Labr}).


\section{Dehornoy's ordering}

Call a group or a monoid $G$ {\it left-orderable} if there exists a strict linear ordering $<$ of 
the elements of $G$ which is left invariant, namely, for $x,y,z \in G$, $x<y$ implies $zx<zy$. 
If, in addition, the ordering is a well-ordering, then we say that $G$ is {\it left-well-orderable}. 
If there is an ordering which is invariant under multiplication on both sides, we say 
that $G$ is {\it biorderable}.

\bigskip\noindent
Recall that the braid group $B_n$ has a presentation with generators $\sigma_1, \dots, \sigma_{n-
1}$ and relations
\begin{gather}
\label{eq41} \sigma_i \sigma_j = \sigma_j \sigma_i \quad \text{for }|i-j|>1\,,\\
\label{eq42} \sigma_i \sigma_j \sigma_i = \sigma_j \sigma_i \sigma_j \quad \text{for }|i-j|=1\,.
\end{gather}
The submonoid generated by $\sigma_1, \dots, \sigma_{n-1}$ is called the {\it positive braid 
monoid} on $n$ strings, and is denoted by $B_n^+$. It has a monoid presentation with generators 
$\sigma_1, \dots, \sigma_{n-1}$ and relations \eqref{eq41} and \eqref{eq42}, it is a Garside 
monoid, and the group of fractions of $B_n^+$ is $B_n$.

\bigskip\noindent
The starting point of the present section is the following result due to Dehornoy \cite{Deh3}.

\bigskip\noindent
{\bf Theorem 3.1 (Dehornoy \cite{Deh3}).} {\it The Artin braid group $B_n$ is left-orderable, by 
an ordering which is a well-ordering when restricted to $B_n^+$.}

\bigskip\noindent
Left-invariant orderings are important because of a long-standing conjecture in group theory: 
that the group ring of a torsion-free group has no zero divisors. This conjecture is true for 
left-orderable groups, thus Theorem 3.1 implies that $\Z [B_n]$ has no zero divisors. However, the 
main point of interest of Theorem 3.1 is not only the mere existence of orderings on the braid 
groups, but the particular nature of the construction. This has been already used to prove that 
some representations of the braid groups by automorphisms of groups are faithful (see \cite{Shp} 
and \cite{CrPa2}), and to study the so-called ``weak faithfulness'' of the Burau representation 
(see \cite{Deh4}), and I am sure that these ideas will be exploited again in the future.

\bigskip\noindent
Let $G$ be a left-orderable group, and let $<$ be a left invariant linear ordering on $G$. Call 
an element $g \in G$ {\it positive} if $g>1$, and denote by $\PP$ the set of positive elements. 
Then $\PP$ is a subsemigroup ({\it i.e.} $\PP \cdot \PP \subset \PP$), and we have the disjoint 
union $G= \PP \sqcup \PP^{-1} \sqcup \{1\}$, where $\PP^{-1}= \{g^{-1}; g \in \PP\}$. Conversely, 
if $\PP$ is a subsemigroup of $G$ such that $G=\PP \sqcup \PP^{-1} \sqcup \{1\}$, then the 
relation $<$ defined by $f<g$ if $f^{-1}g \in \PP$ is a left invariant linear ordering on $G$.

\bigskip\noindent
The description of the positive elements in Dehornoy's ordering is based on the following result 
which, in my view, is more interesting than the mere existence of the ordering, 
because it leads to new techniques to study the braid groups. 

\bigskip\noindent
Let $B_{n-1}$ be the subgroup of $B_n$ generated by $\sigma_1, \dots, \sigma_{n-2}$. Let $\beta 
\in B_n$. We say that $\beta$ is {\it $\sigma_{n-1}$-positive} if it can be written
\[
\beta= \alpha_0 \sigma_{n-1} \alpha_1 \sigma_{n-1} \alpha_2 \dots \sigma_{n-1} \alpha_l \,,
\]
where $l \ge 1$ and $\alpha_0, \alpha_1, \dots, \alpha_l \in B_{n-1}$. We say that $\beta$ 
is {\it $\sigma_{n-1}$-negative} if it can be written
\[
\beta= \alpha_0 \sigma_{n-1}^{-1} \alpha_1 \sigma_{n-1}^{-1} \alpha_2 \dots \sigma_{n-1}^{-1} 
\alpha_l\,,
\]
where $l \ge 1$ and $\alpha_0, \alpha_1, \dots, \alpha_l \in B_{n-1}$. We denote by $\PP_{n-1}$ 
the set of $\sigma_{n-1}$-positive elements, and by $\PP_{n-1}^{-}$ the set of 
$\sigma_{n-1}$-negative elements. Note that $\PP_{n-1}^{-}=\PP_{n-1}^{-1}=\{ \beta^{-1}; \beta \in \PP_{n-1}\}$.

\bigskip\noindent
{\bf Theorem 3.2 (Dehornoy \cite{Deh3}).} {\it We have the disjoint union $B_n= \PP_{n-1} \sqcup 
\PP_{n-1}^{-} \sqcup B_{n-1}$.}

\bigskip\noindent
Now, the set $\PP$ of positive elements in Dehornoy's ordering can be described as follows. For 
$k=1, \dots, n-1$, let $\PP_k$ denote the set of $\sigma_k$-positive elements in $B_{k+1}$, 
where $B_{k+1}$ is viewed as the subgroup of $B_n$ generated by $\sigma_1, \dots, \sigma_k$. Then 
$\PP= \PP_1 \sqcup \PP_2 \sqcup \dots \sqcup \PP_{n-1}$. It is easily checked using Theorem 3.2 
that $\PP$ is a subsemigroup of $B_n$ and that $B_n= \PP \sqcup \PP^{-1} \sqcup \{1\}$.

\bigskip\noindent
Dehornoy's proof of Theorem 3.2 involves rather delicate combinatorial and algebraic 
constructions which were partly motivated by questions in set theory. A topological construction 
of Dehornoy's ordering of the braid group $B_n$, viewed as the mapping class group of the 
$n$-punctured disk, is given in \cite{FGRRW}. This last construction has been extended to the mapping 
class groups of the punctured surfaces with non-empty boundary by Rourke and Wiest \cite{RoWi}. 
This extension involves some choices, and, therefore, is not canonical. However, it is strongly 
connected to Mosher's normal forms \cite{Mos} and, as with Dehornoy's ordering, is sufficiently 
explicit to be used for other purposes (to study faithfulness properties of ``some'' 
representations of the mapping class group, for example. This is just an idea...). Another 
approach (see \cite{ShWi}) based on Nielsen's ideas \cite{Nie} gives rise to many other linear 
orderings of the braid groups, and, more generally, of the mapping class groups of punctured 
surfaces with non-empty boundary, but I am not convinced that these orderings are explicit enough 
to be used for other purposes.

\bigskip\noindent
The mapping class group of a (punctured) closed surface is definitively not left-orderable 
because it has torsion. However, it is an open question to determine whether some special 
torsion-free subgroups of $\MM( \Sigma)$ such as Ivanov's group 
$\Gamma_m(\Sigma)$ (see \cite{Iva3}) are orderable. 
The Torelli group is residually torsion-free nilpotent
(see \cite{Hain}), therefore is biorderable.
We refer to \cite{DDRW} for a detailed 
account on Dehornoy's ordering and related questions.


\section{Thurston classification}

Recall that an {\it essential circle} in $\Sigma \setminus \PP$ is an oriented embedding $a: \S^1 
\to \Sigma \setminus \PP$ such that $a(\S^1) \cap \partial \Sigma = \emptyset$, and $a(\S^1)$ 
does not bound a disk in $\Sigma \setminus \PP$. An essential circle $a: \S^1 \to \Sigma 
\setminus \PP$ is called {\it generic} if it does not bound a disk in $\Sigma$ containing one 
puncture, and if it is not isotopic to a boundary component of $\Sigma$. Recall that two circles 
$a,b: \S^1 \to \Sigma \setminus \PP$ are {\it isotopic} if there exists a continuous family $a_t: 
\S^1 \to \Sigma \setminus \PP$, $t \in [0,1]$, such that $a_0=a$ and $a_1=b$. Isotopy of circles 
is an equivalence relation that we denote by $a \sim b$. Observe that $h(a) \sim h'(a')$ if $a 
\sim a'$ and $h,h' \in \HH( \Sigma, \PP)$ are isotopic. So, the mapping class group 
$\MM(\Sigma,\PP)$ acts on the set $\CC (\Sigma, \PP)$ of isotopy classes of generic circles.

\bigskip\noindent
Let $f \in \MM( \Sigma, \PP)$. We say that $f$ is a {\it pseudo-Anosov element} if it has no 
finite orbit in $\CC(\Sigma, \PP)$, we say that $f$ is {\it periodic} if $f^m$ acts trivially on 
$\CC (\Sigma, \PP)$ for some $m \ge 1$, and we say that $f$ is {\it reducible} otherwise.

\bigskip\noindent
The use of the action of $\MM(\Sigma, \PP)$ on $\CC (\Sigma, \PP)$ and the ``Thurston 
classification'' of the elements of $\MM(\Sigma, \PP)$ play a prominent role in the study of the 
algebraic properties of the mapping class groups. Here are some applications.

\bigskip\noindent
{\bf Theorem 4.1.}
\begin{enumerate}
\item {\bf (Ivanov \cite{Iva3}).} {\it Let $f$ be a pseudo-Anosov element of $\MM(\Sigma, \PP)$. 
Then $f$ has infinite order, $\langle f \rangle \cap Z( \MM(\Sigma, \PP)) = \{1\}$, where 
$\langle f \rangle$ denotes the subgroup generated by $f$, and $\langle f \rangle \times 
Z(\MM(\Sigma,\PP))$ is a subgroup of finite index of the centralizer of $f$ in $\MM(\Sigma, 
\PP)$.}
\item ({\bf Birman, Lubotzky, McCarthy \cite{BLM}}). {\it If a subgroup of $\MM(\Sigma, \PP)$ 
contains a solvable subgroup of finite index, then it contains an abelian subgroup of finite 
index.}
\item ({\bf Ivanov \cite{Iva3}, McCarthy \cite{McM}}). {\it Let $G$ be a subgroup of $\MM(\Sigma, 
\PP)$. Then either $G$ contains a free group of rank 2, or $G$ contains an abelian subgroup of 
finite index.}
\end{enumerate}

\bigskip\noindent
The use of the action of $\MM(\Sigma, \PP)$ on $\CC (\Sigma, \PP)$ is also essential in the 
determination of the (maximal) free abelian subgroups of $\MM(\Sigma, \PP)$ (see \cite{BLM}), and 
in the calculation of the automorphism group and abstract commensurator group of $\MM(\Sigma, \PP)$ 
(see \cite{Iva1}, \cite{Iva2}, \cite{McM2}, \cite{IvMc}).

\bigskip\noindent
The present section concerns the question of how to translate these techniques to other groups 
such as the (spherical type) Artin groups, or the Coxeter groups.

\bigskip\noindent
Let $\Gamma$ be a Coxeter graph, let $M=(m_{s\,t})_{s,t\in S}$ be the Coxeter matrix of $\Gamma$,
and let $A=A_\Gamma$ be the Artin group of type $\Gamma$. For $X \subset S$, we denote by $\Gamma_X$
the full subgraph of $\Gamma$ generated by $X$, and by $A_X$ the subgroup of $A$ generated by
$\SS_X= \{ \sigma_x; x \in X\}$.
By \cite{Lek} (see also \cite{Par3}), the 
group $A_X$ is the Artin group of type $\Gamma_X$. A subgroup of the form $A_X$ is called a {\it 
standard parabolic subgroup}, and a subgroup which is conjugate to a standard parabolic subgroup is simply 
called a {\it parabolic subgroup}. We shall denote by $\AA\CC=\AA\CC(\Gamma)$ the set of 
parabolic subgroups different from $A$ and from $\{1\}$.

\bigskip\noindent
Let $f \in A$. We say that $f$ is an {\it (algebraic) pseudo-Anosov element} if it has no finite 
orbit in $\AA\CC$, we say that $f$ is an {\it (algebraic) periodic element} if $f^m$ acts 
trivially on $\AA\CC$ for some $m >0$, and we say that $f$ is an {\it (algebraic) reducible 
element} otherwise.

\bigskip\noindent
The above definitions are motivated by the following result.

\bigskip\noindent
{\bf Proposition 4.2.} {\it Let $f$ be an element of the braid group $B_n$ (viewed as the mapping 
class group $\MM(\D, \PP)$ as well as the Artin group associated to the Dynkin graph $A_{n-
1}$).
\begin{itemize}
\item $f$ is a pseudo-Anosov element if and only if $f$ is an algebraic pseudo-Anosov element.
\item $f$ is a periodic element if and only if $f$ is an algebraic periodic element.
\item $f$ is a reducible element if and only if $f$ is an algebraic reducible element.
\end{itemize}}

\bigskip\noindent
{\bf Proof.} Recall that $\D= \{ z \in \C; |z| \le 1\}$, $\PP=\{P_1, \dots, P_n\}$ is a set of 
$n$ punctures in the interior of $\D$, and $B_n= \MM(\D, \PP)$. Let $\sigma_1, \dots, \sigma_{n-
1}$ be the standard generators of $B_n$. We can assume that each $\sigma_i$ acts trivially 
outside a disk containing $P_i$ and $P_{i+1}$, and permutes $P_i$ and $P_{i+1}$.

\bigskip\noindent
Let $a: \S^1 \to \D \setminus \PP$ be a generic circle. Then $a(\S^1)$ separates $\D$ into two 
components: a disk $\D_a$ bounded by $a(\S^1)$, and an annulus $\A_a$ bounded by $a(\S^1) \sqcup 
\S^1$, where $\S^1= \{z \in \C; |z|=1\}$ (see Figure 4). Let $m= |\D_a \cap \PP|$. Then the 
hypothesis ``$a$ is generic'' implies that $2 \le m\le n-1$. Furthermore, the group $\MM(\D_a, 
\D_a \cap \PP)$ is isomorphic to $B_m$.

\begin{figure}
\begin{center}
\includegraphics[width=5cm]{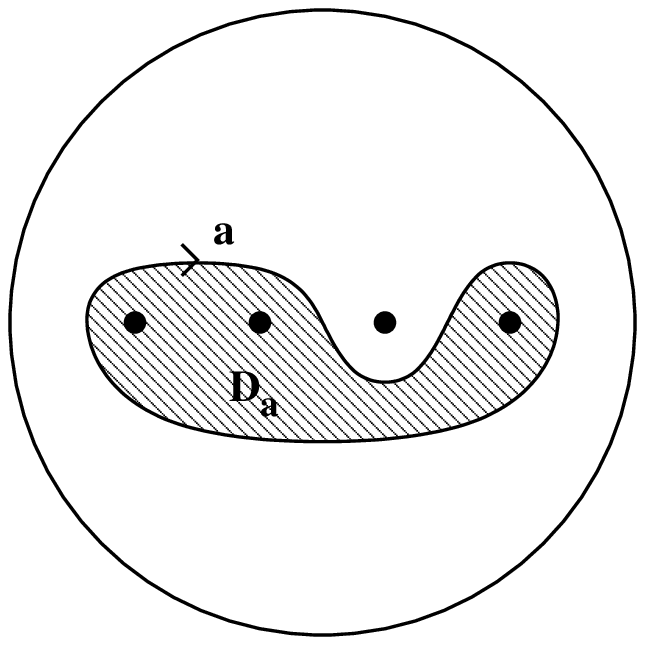}
\end{center}
\caption{The disk $\D_a$.}
\end{figure}

\bigskip\noindent
Let $\tau_a$ denote the Dehn twist along $a$. If $m=|\D_a \cap \PP| \ge 3$, then the center $Z(\MM( \D_a, \D_a 
\cap \PP))$ of $\MM( \D_a, \D_a \cap \PP)$ is the (infinite) cyclic subgroup generated by 
$\tau_a$. If $m=2$, then $\MM( \D_a, \D_a \cap \PP)$ is an infinite cyclic group, and $\tau_a$ 
generates the (unique) index 2 subgroup of $\MM (\D_a, \D_a \cap \PP)$. Recall that, if 
$\tau_{a_1}$ and $\tau_{a_2}$ are two Dehn twists along generic circles $a_1$ and $a_2$, 
respectively, and if $\tau_{a_1}^{k_1} = \tau_{a_2}^{k_2}$ for some $k_1,k_2 \in \Z \setminus \{0\}$, 
then $k_1=k_2$ and $a_1=a_2$. Recall also that, if $\tau_a$ is the Dehn twist along a generic circle $a$, then 
$g \tau_a g^{-1}$ is the Dehn twist along $g(a)$, for any $g \in \MM(\D, \PP)$. We conclude from 
the above observations that, if $g \cdot \MM( \D_a, \D_a \cap \PP) \cdot g^{-1} = \MM (\D_a, 
\D_a \cap \PP)$, then $g(a)=a$. Conversely, if $g(a)=a$, then $g(\D_a, \D_a \cap \PP)= (\D_a, 
\D_a \cap \PP)$, thus $g \cdot \MM( \D_a, \D_a \cap \PP) \cdot g^{-1} = \MM (\D_a, \D_a \cap 
\PP)$.

\bigskip\noindent
Write $S=\{1,2, \dots, n-1\}$. For $X \subset S$, we set $\SS_X= \{\sigma_x; x\in X\}$, and 
denote by $(B_n)_X$ the subgroup of $B_n$ generated by $\SS_X$.

\bigskip\noindent
Let $1\le m\le n-2$, and $X(m)=\{ 1,2, \dots, m\}$. We can choose a generic circle $a_m: \S^1 
\to \D \setminus \PP$ such that $\D_{a_m} \cap \PP = \{P_1, \dots, P_m, P_{m+1}\}$ and $\MM( 
\D_{a_m}, \D_{a_m} \cap \PP)= (B_n)_{X(m)}$ (see Figure 5). Observe that, if $a: \S^1 \to \D 
\setminus \PP$ is a generic circle such that $|\D_a \cap \PP|=m+1$, then there exists $g \in \MM( 
\D, \PP)= B_n$ such that $g \MM (\D_a, \D_a \cap \PP) g^{-1} = \MM (\D_{a_m}, \D_{a_m} \cap \PP)= 
(B_n)_{X(m)}$. In particular, $\MM (\D_a, \D_a \cap \PP)$ is a parabolic subgroup of $B_n$.

\begin{figure}
\begin{center}
\includegraphics[width=5cm]{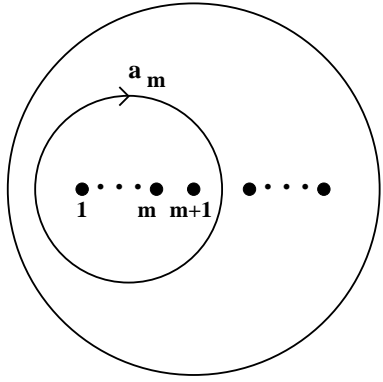}
\end{center}
\caption{The curve $a_m$.}
\end{figure}

\bigskip\noindent
Take $X \subset S$, $X \neq S$. Decompose $X$ as a disjoint union $X=X_1 \sqcup X_2 \sqcup \dots 
\sqcup X_q$, where $X_k= \{ i_k, i_k+1, \dots, i_k+t_k\}$ for some $i_k \in \{1, \dots, n-1\}$ 
and some $t_k \in \N$, and $i_k +t_k +1< i_{k+1}$, for all $k=1, \dots, q$. We can choose generic 
circles $b_1, \dots, b_q: \S^1 \to \D \setminus \PP$ such that
\begin{gather*}
\D_{b_k} \cap \PP = \{P_{i_k}, \dots, P_{i_k+t_k}, P_{i_k+t_k+1} \}\,, \quad \text{for } 1\le 
k\le q\,,\\
\D_{b_k} \cap \D_{b_l} = \emptyset\,, \quad \text{for }1\le k \neq l\le q\,,\\
\MM( \D_{b_k}, \D_{b_k} \cap \PP)= (B_n)_{X_k} \quad \text{for } 1 \le k \le q\,.
\end{gather*}
(See Figure 6.) It is easily checked that
\[
(B_n)_X= (B_n)_{X_1} \times \dots \times (B_n)_{X_q}= \MM(\D_{b_1}, \D_{b_1} \cap \PP) \times 
\dots \times \MM (\D_{b_q}, \D_{b_q} \cap \PP)\,,
\]
and the only Dehn twists which lie in the center of $(B_n)_X$ are $\tau_{b_1}, \tau_{b_2}, \dots, 
\tau_{b_q}$.

\begin{figure}
\begin{center}
\includegraphics[width=5cm]{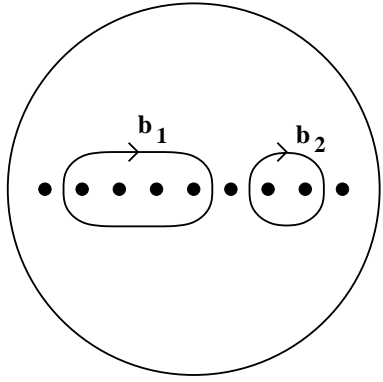}
\end{center}
\caption{The curves $b_1$, $b_2$ for $X=\{2,3,4,7\}$.}
\end{figure}

\bigskip\noindent
Let $g \in \MM(\D, \PP)$ which fixes some $a \in \CC (\D, \PP)$. Then $g \MM( \D_a, \D_a \cap 
\PP) g^{-1} = \MM (\D_a, \D_a \cap \PP)$ and $\MM (\D_a, \D_a \cap \PP)$ is a parabolic subgroup, 
thus $g$ fixes some element in $\AA\CC$.

\bigskip\noindent
Let $g \in \MM (\D, \PP)$ which fixes some element in $\AA\CC$. Up to conjugation, we may suppose 
that there exists $X \subset S$ such that $g (B_n)_X g^{-1} = (B_n)_X$. We take again the 
notations introduced above. Recall that the only Dehn twists which lie in the center of $(B_n)_X$ 
are $\tau_{b_1}, \tau_{b_2}, \dots, \tau_{b_q}$, thus $g$ must permute the set $\{ \tau_{b_1}, 
\tau_{b_2}, \dots, \tau_{b_q}\}$. So, there exists $m\ge 1$ such that $g^m \tau_{b_1} g^{-m} = 
\tau_{b_1}$, that is, $g^m(b_1)=b_1$.

\bigskip\noindent
Let $g \in \MM( \D, \PP)$ which fixes all the elements of $\CC(\D, \PP)$. Let $A \in \AA\CC$. 
There exist $h \in B_n$ and $X \subset S$ such that $A=h (B_n)_X h^{-1}$. We keep the notations 
used above. The element $h^{-1}gh$ also fixes all the elements of $\CC (\D, \PP)$. In particular, 
we have $(h^{-1}gh) (b_k)=b_k$ for all $k =1, \dots, q$. This implies that $h^{-1}gh$ fixes $\MM 
(\D_{b_k}, \D_{b_k} \cap \PP)$ for all $k=1, \dots, q$, thus $h^{-1}gh$ fixes $(B_n)_X= \MM( 
\D_{b_1}, \D_{b_1} \cap \PP) \times \dots \times \MM (\D_{b_q}, \D_{b_q} \cap \PP)$, therefore 
$g$ fixes $A=h (B_n)_X h^{-1}$.

\bigskip\noindent
Let $g \in \MM (\D, \PP)$ which fixes all the elements of $\AA\CC$. Let $a \in \CC( \D, \PP)$. 
Since $\MM( \D_a, \D_a \cap \PP)$ is a parabolic subgroup, we conclude that $g \MM (\D_a, \D_a 
\cap \PP) g^{-1} = \MM (D_a, \D_a \cap \PP)$, thus $g(a)=a$.
\endproof

\bigskip\noindent
In spite of the fact that the normalizers of the parabolic subgroups of the spherical type Artin 
groups are quite well understood (see \cite{Par3}, \cite{God}), it is not known, for instance, 
whether the centralizer of an algebraic pseudo-Anosov element $g$ in a spherical type Artin group 
$A$ has $\langle g \rangle \times Z(A)$ as finite index subgroup. This question might be the 
starting point for understanding the abelian subgroups of the spherical type Artin groups.

\bigskip\noindent
We turn now to a part of Krammer's Ph. D. thesis which describes
a phenomenom for Coxeter groups similar to the ``Thurston
classification'' for the mapping class groups: in a non-affine irreducible
Coxeter group $W$, there are certain elements, called {\it essential}
elements, which have infinite order, and have the property that
the cyclic subgroups generated by them have finite index in
their centralizer. Although almost all elements of $W$ are
essential, it is very hard to determine whether a given element is essential.

\bigskip\noindent
Let $\Gamma$ be a Coxeter graph with set of vertices $S$, and let $W=W_\Gamma$ be the
Coxeter group of type $\Gamma$.
Let $\Pi= \{ \alpha_s; s \in S\}$ be an abstract set in one-to-one correspondence with $S$. The 
elements of $\Pi$ are called {\it simple roots}. Let $U= \oplus_{s \in S} \R \alpha_s$ denote the 
(abstract) real vector space having $\Pi$ as a basis. Define the {\it canonical form} as the 
symmetric bilinear form $\langle, \rangle: U \times U \to \R$ determined by
\[
\langle \alpha_s, \alpha_t \rangle = \left\{
\begin{array}{ll}
-\cos (\pi/m_{s\,t}) \quad &\text{if }m_{s\,t}<+\infty\,,\\
-1\quad&\text{if } m_{s\,t}=+\infty\,.
\end{array}\right.
\]
For $s \in S$, define the linear transformation $\rho_s: U \to U$ by
\[
\rho_s(x)=x-2\langle x,\alpha_s \rangle \alpha_s\,, \quad \text{for }x \in U\,.
\]
Then the map $s \mapsto \rho_s$, $s \in S$, determines a well-defined representation $W \to 
GL(U)$, called {\it canonical representation}. This representation is faithful and preserves the 
canonical form.

\bigskip\noindent
The set $\Phi= \{ w \alpha_s; w \in W \text{ and } s \in S\}$ is called the {\it root system} of 
$W$. The set $\Phi^+= \{ \beta= \sum_{s \in S} \lambda_s \alpha_s \in \Phi; \lambda_s \ge 0\text{ 
for all } s \in S\}$ is the set of {\it positive roots}, and $\Phi^-= -\Phi^+$ is the set of {\it 
negative roots}. For $w \in W$, we set
\[
\Phi_w= \{ \beta \in \Phi^+\ ;\ w\beta \in \Phi^-\}\,.
\]
The following proposition is a mixture of several well-known facts on $\Phi$. The proofs 
can be found in \cite{Hil}.

\bigskip\noindent
{\bf Proposition 4.3.}{\it
\begin{enumerate}
\item We have the disjoint union $\Phi= \Phi^+ \sqcup \Phi^-$.
\item Let $\lg: W \to\N$ denote the word length of $W$ with respect to $S$. Then $\lg(w)= |\Phi_w|$ 
for all $w \in W$. In particular, $\Phi_w$ is finite.
\item Let $w \in W$ and $s \in S$. Then
\[
\lg(ws)=\left\{
\begin{array}{ll}
\lg(w)+1\quad&\text{if } w\alpha_s \in \Phi^+\,,\\
\lg(w)-1\quad&\text{if }w\alpha_s \in \Phi^-\,.
\end{array}\right.
\]
\item Let $w \in W$ and $s \in S$. Write $\beta= w \alpha_s \in \Phi$, and $r_\beta = w s w^{-1} 
\in W$. Then $r_\beta$ acts on $U$ by
\[
r_\beta (x)= x-2 \langle x, \beta \rangle \beta\,.
\]
\end{enumerate}}

\noindent
Let $u,v \in W$ and $\alpha \in \Phi$. We say that $\alpha$ {\it separates} $u$ and $v$ if there exists 
$\varepsilon \in \{ \pm 1\}$ such that $u \alpha \in \Phi^\varepsilon$ and $v \alpha \in \Phi^{-
\varepsilon}$. Let $w \in W$ and $\alpha \in \Phi$. We say that $\alpha$ is {\it $w$-periodic} if 
there exists some $m \ge 1$ such that $w^m \alpha = \alpha$.

\bigskip\noindent
{\bf Proposition 4.4.} {\it Let $w \in W$ and $\alpha \in \Phi$. Then precisely one of the 
following holds.
\begin{enumerate}
\item $\alpha$ is $w$-periodic.
\item $\alpha$ is not $w$-periodic, and the set $\{m \in \Z; \alpha \text{ separates } w^m \text{ 
and } w^{m+1} \}$ is finite and even.
\item $\alpha$ is not $w$-periodic, and the set $\{m \in \Z; \alpha \text{ separates } w^m \text{ 
and } w^{m+1} \}$ is finite and odd.
\end{enumerate}}

\bigskip\noindent
Call $\alpha$ {\it $w$-even} in Case 2, and {\it $w$-odd} in Case 3.

\bigskip\noindent
{\bf Proof.} Assume that $\alpha$ is not $w$-periodic, and put $N_w(\alpha)= \{m \in \Z; \alpha \text{ 
separates } w^m \text{ and } w^{m+1} \}$. We have to show that $N_w(\alpha)$ is finite.

\bigskip\noindent
If $m \in N_w(\alpha)$, then $w^m \alpha \in \Phi_w \cup -\Phi_w$. On the other hand, if $w^{m_1} 
\alpha = w^{m_2} \alpha$, then $w^{m_1-m_2} \alpha = \alpha$, thus $m_1-m_2=0$, since $\alpha$ is 
not $w$-periodic. Since $\Phi_w \cup -\Phi_w$ is finite (see Proposition 4.3), we conclude that 
$N_w(\alpha)$ is finite.
\endproof

\bigskip\noindent
For $X \subset S$, we denote by $W_X$ the subgroup of $W$ generated by $X$. Such a subgroup is 
called a {\it standard parabolic subgroup}. A subgroup of $W$ conjugated to a standard parabolic 
subgroup is simply called a {\it parabolic subgroup}. An element $w \in W$ is called {\it 
essential} if it does not lie in any proper parabolic subgroup. Finally, recall that a Coxeter 
group $W$ is said to be of {\it affine type} if the canonical form $\langle, \rangle$ is 
non-negative (namely, $\langle x,x \rangle \ge 0$ for all $x \in U$).

\bigskip\noindent
Now, we can state Krammer's result.

\bigskip\noindent
{\bf Theorem 4.5 (Krammer \cite{Kra}).} {\it Assume $W$ to be an irreducible non-affine Coxeter 
group. Let $w \in W$.
\begin{enumerate}
\item $w$ is essential if and only if $W$ is generated by the set $\{ r_\beta; \beta \in \Phi^+ 
\text{ and } \beta\ w\text{-odd}\}$.
\item Let $m \in \N$, $m \ge 1$. $w$ is essential if and only if $w^m$ is essential.
\item Suppose $w$ is essential. Then $w$ has infinite order, and $\langle w \rangle = \{ w^m; m 
\in \Z\}$ is a finite index subgroup of the centralizer of $w$ in $W$.
\end{enumerate}}



\bigskip\bigskip\noindent
{\bf Luis Paris},

\smallskip\noindent 
Institut de Math\'ematiques de Bourgogne, UMR 5584 du CNRS, Universit\'e de Bourgogne, B.P. 
47870, 21078 Dijon cedex, France

\smallskip\noindent
E-mail: {\tt lparis@u-bourgogne.fr}


\begin{thebibliography}{99}

\bibitem{Arn}
V.I. Arnold,
{\it Normal forms of functions near degenerate critical points, the Weyl groups $A_k$, $D_k$, 
$E_k$ and Lagrangian singularities},
Funkcional. Anal. i Prilov zen. {\bf 6} (1972), no. 4, 3--25.

\bibitem{BesCor}
D. Bessis, R. Corran,
{\it Garside structure for the braid group of $G(e,e,r)$},
preprint.

\bibitem{Bir1}
J.S. Birman,
{\it Braids, links, and mapping class groups},
Annals of Mathematics Studies, No. 82, Princeton University Press,
Princeton, N.J., 1974.

\bibitem{Bir2}
J.S. Birman,
{\it Mapping class groups of surfaces},
Braids (Santa Cruz, CA, 1986), 13--43, Contemp. Math., 78, Amer. Math.
Soc., Providence, RI, 1988.

\bibitem{BLM}
J.S. Birman, A. Lubotzky, J. McCarthy,
{\it Abelian and solvable subgroups of the mapping class groups},
Duke Math. J. {\bf 50} (1983), 1107--1120.

\bibitem{Bou}
N. Bourbaki,
{\it Groupes et alg\`ebres de Lie, chapitres 4, 5 et 6},
Hermann, Paris, 1968.

\bibitem{Bri}
E. Brieskorn,
{\it Die Fundamentalgruppe des Raumes der regul\"aren Orbits einer endlichen komplexen 
Spiegelungsgruppe},
Invent. Math. {\bf 12} (1971), 57--61.

\bibitem{BrSa}
E. Brieskorn, K. Saito,
{\it Artin-Gruppen und Coxeter-Gruppen},
Invent. Math. {\bf 17} (1972), 245--271.

\bibitem{ChDa}
R. Charney, M.W. Davis,
{\it The $K(\pi,1)$-problem for hyperplane complements associated to infinite reflection groups},
J. Amer. Math. Soc. {\bf 8} (1995), 597--627.

\bibitem{CMW}
R. Charney, J. Meier, K. Whittlesey,
{\it Bestvina's normal form complex and the homology of Garside groups},
Geom. Dedicata {\bf 105} (2004), 171--188.

\bibitem{CrPa}
J. Crisp, L. Paris,
{\it The solution to a conjecture of Tits on the subgroup generated by the squares of the 
generators of an Artrin group},
Invent. Math. {\bf 145} (2001), 19--36.

\bibitem{CrPa2}
J. Crisp, L. Paris,
{\it Representations of the braid group by automorphisms of groups, invariants of links, and 
Garside groups},
Pacific J. Math., to appear.

\bibitem{Deh3}
P. Dehornoy,
{\it Braid groups and left distributive operations},
Trans. Amer. Math. Soc. {\bf 343} (1994), 115--150.

\bibitem{Deh4}
P. Dehornoy,
{\it Weak faithfulness properties for the Burau representation},
Topology Appl. {\bf 69} (1996), 121--143.

\bibitem{Deh2}
P. Dehornoy,
{\it Gaussian groups are torsion free},
J. Algebra {\bf 210} (1998), 291--297.

\bibitem{Deh1}
P. Dehornoy,
{\it Groupes de Garside},
Ann. Sci. Ecole Norm. Sup. (4) {\bf 35} (2002), 267--306.

\bibitem{DDRW}
P. Dehornoy, I. Dynnikov, D. Rolfsen, B. Wiest,
{\it Why are braids orderable?}
Panoramas et Synth\`eses, 14, Soci\'et\'e Math\'ematique de France, Paris, 2002.

\bibitem{DeLa}
P. Dehornoy, Y. Lafont,
{\it Homology of Gaussian groups},
Ann. Inst. Fourier (Grenoble) {\bf 53} (2003), 489--540.

\bibitem{DePa}
P. Dehornoy, L. Paris,
{\it Gaussian groups and Garside groups, two generalisations of Artin groups},
Proc. London Math. Soc. (3) {\bf 79} (1999), 569--604.

\bibitem{Del}
P. Deligne,
{\it Les immeubles des groupes de tresses g\'en\'eralis\'es},
Invent. Math. {\bf 17} (1972), 273--302.

\bibitem{DuKr}
G. Duchamp, D. Krob,
{\it The lower central series of the free partially commutative group},
Semigroup Forum {\bf 45} (1992), 385--394.

\bibitem{EpAl}
D.B.A. Epstein, J.W. Cannon, D.F. Holt, S.V.F. Levy, M.S. Paterson, W.P. Thurston,
{\it Word processing in groups},
Jones and Bartlett Publishers, Boston, MA, 1992.

\bibitem{FaNe}
E. Fadell, L. Neuwirth,
{\it Configuration spaces},
Math. Scand. {\bf 10} (1962), 111--118.

\bibitem{FGRRW}
R. Fenn, M.T. Greene, D. Rolfsen, C. Rourke, B. Wiest,
{\it Ordering the braid groups},
Pacific J. Math. {\bf 191} (1999), 49--74.

\bibitem{FrGo}
N. Franco, J. Gonz\'alez-Meneses,
{\it Conjugacy problem for braid groups and Garside groups},
J. Algebra {\bf 266} (2003), 112--132.

\bibitem{Gar}
F.A. Garside,
{\it The braid group and other groups},
Quart. J. Math. Oxford Ser. (2) {\bf 20} (1969), 235--254.

\bibitem{Hain}
R. Hain,
{\it Infinitesimal presentations of the Torelli groups},
J. Amer. Math. Soc. {\bf 10} (1997), no. 3, 597--651.

\bibitem{God}
E. Godelle,
{\it Normalisateur et groupe d'Artin de type sph\'erique},
J. Algebra {\bf 269} (2003), 263--274.

\bibitem{Hil}
H. Hiller,
{\it Geometry of Coxeter groups},
Research Notes in Mathematics, 54, Pitman (Advanced Publishing Program), Boston, 
Mass.--London, 1982.

\bibitem{Hum}
J.E. Humphreys,
{\it Reflection groups and Coxeter groups},
Cambridge Studies in Advanced Mathematics, 29,
Cambridge University Press, Cambridge, 1990.

\bibitem{Iva1}
N.V. Ivanov,
{\it Automorphisms of Teichm\"uller modular groups},
Topology and geometry - Rohlin Seminar, 
199--270, Lecture Notes in Math., 1346,
Springer, Berlin, 1988.

\bibitem{Iva3}
N.V. Ivanov,
{\it Subgroups of Teichm\"uller modular groups},
Translations of Mathematical Monographs, 115, American Mathematical Society, Providence, RI, 
1992.

\bibitem{Iva2}
N.V. Ivanov,
{\it Automorphism of complexes of curves and of Teichm\"uller spaces},
Internat. Math. Res. Notices {\bf 14} (1997), 651--666.

\bibitem{IvMc}
N.V. Ivanov, J.D. McCarthy,
{\it On injective homomorphisms between Teichm\"uller modular groups. I},
Invent. Math. {\bf 135} (1999), 425--486.

\bibitem{Kra}
D. Krammer,
{\it The conjugacy problem for Coxeter groups},
Ph. D. Thesis, Universiteit Utrecht, 1995.

\bibitem{Labr}
C. Labru\`ere,
{\it Generalized braid groups and mapping class groups},
J. Knot Theory Ramifications {\bf 6} (1997), 715--726.

\bibitem{LaPa}
C. Labru\`ere, L. Paris,
{\it Presentations for the punctured mapping class groups in terms of Artin groups},
Algebr. Geom. Topol. {\bf 1} (2001), 73--114.

\bibitem{Lek}
H. Van der Lek,
{\it The homotopy type of complex hyperplane complements},
Ph. D. Thesis, Nijmegen, 1983.

\bibitem{Mat}
M. Matsumoto,
{\it A presentation of mapping class groups in terms of Artin groups and geometric monodromy of 
singularities},
Math. Ann. {\bf 316} (2000), 401--418.

\bibitem{McM}
J. McCarthy,
{\it A ``Tits-alternative'' for subgroups of surface mapping class groups},
Trans. Amer. Math. Soc. {\bf 291} (1985), 583--612.

\bibitem{McM2}
J. McCarthy,
{\it Automorphisms of surface mapping class groups. A recent theorem of N. Ivanov},
Invent. Math. {\bf 84} (1986), 49--71.

\bibitem{Mos}
L. Mosher,
{\it Mapping class groups are automatic},
Ann. of Math. (2) {\bf 142} (1995), 303--384.

\bibitem{Nie}
J. Nielsen,
{\it Untersuchungen zur Topologie des geschlossenen zweiseitigen Fl\"achen},
Acta Math. {\bf 50} (1927), 189--358.

\bibitem{Par3}
L. Paris,
{\it Parabolic subgroups of Artin groups},
J. Algebra {\bf 196} (1997), 369--399.

\bibitem{PeVa}
B. Perron, J.P. Vannier,
{\it Groupe de monodromie g\'eom\'etrique des singularit\'es simples},
Math. Ann. {\bf 306} (1996), 231--245.

\bibitem{Pic2}
M. Picantin,
{\it The conjugacy problem in small Gaussian groups},
Comm. Algebra {\bf 29} (2001), 1021--1039.

\bibitem{Pic1}
M. Picantin,
{\it Automatic structures for torus link groups},
J. Knot Theory Ramifications {\bf 12} (2003), 833--866.

\bibitem{RoWi}
C. Rourke, B. Wiest,
{\it Order automatic mapping class groups},
Pacific J. Math. {\bf 194} (2000), 209--227.

\bibitem{ShWi}
H. Short, B. Wiest,
{\it Orderings of mapping class groups after Thurston},
Enseign. Math. (2) {\bf 46} (2000), 279--312.

\bibitem{Shp}
V. Shpilrain,
{\it Representing braids by automorphisms},
Internat. J. Algebra Comput. {\bf 11} (2001), 773--777.

\bibitem{Waj}
B. Wajnryb,
{\it Artin groups and geometric monodromy},
Invent. Math. {\bf 138} (1999), 563--571.

\end{thebibliography}
\end{document}